\magnification=1200

\def\title#1{{\titlefont\noindent #1\bigskip}}

\def\author#1{{\largefont\noindent #1}\medskip}

\def\beginlinemode{\endmode
 \begingroup\obeylines\def\endmode{\par\endgroup}}
\let\endmode=\par

\newbox\theaddress
\def\address{\smallskip\beginlinemode\parindent 0in\getaddress}
{\obeylines
\gdef\getaddress #1 
 #2
 {#1\gdef\addressee{#2}%
   \global\setbox\theaddress=\vbox\bgroup\raggedright%
    \everypar{\hangindent2em}#2
   \def\endaddress{\egroup\endgroup \copy\theaddress \medskip}}}

\def\thanks#1{\footnote{}{\eightpoint #1}}

\long\def\Abstract#1{{\eightpoint\narrower\vskip\baselineskip\noindent
#1\smallskip}}

\def\skipfirstword#1 {}

\def\ir#1{\csname #1\endcsname}

\newdimen\currentht
\newbox\droppedletter
\newdimen\droppedletterwdth
\newdimen\drophtinpts
\newdimen\dropindent

\def\irrnSection#1#2{\edef\tttempcs{\ir{#2}}
\vskip\baselineskip\penalty-3000
{\largefont\bf\noindent \expandafter\skipfirstword\tttempcs. #1}
\vskip6pt}

\def\irSubsection#1#2{\edef\tttempcs{\ir{#2}}
\vskip\baselineskip\penalty-3000
{\bf\noindent \expandafter\skipfirstword\tttempcs. #1}
\vskip6pt}

\def\irSubsubsection#1#2{\edef\tttempcs{\ir{#2}}
\vskip\baselineskip\penalty-3000
{\noindent \expandafter\skipfirstword\tttempcs. #1}
\vskip6pt}

\def\References{\vbox to.25in{\vfil}\noindent{}{\bf References}
\vskip6pt\par}

\def\References{\vskip6pt\noindent{}{\bf References}
\vskip6pt\par}

\def\baselinebreak{\par \ifdim\lastskip<6pt
         \removelastskip\penalty-200\vskip6pt\fi}

\long\def\prclm#1#2#3{\baselinebreak
\noindent{\bf \csname #2\endcsname}:\enspace{\sl #3\par}\baselinebreak}

\def\rem#1#2{\baselinebreak\noindent{\bf \csname #2\endcsname}:}

\def\qed{{$\diamondsuit$}\vskip6pt}

\def\bibitem#1{\par\indent\llap{\rlap{\bf [#1]}\indent}\indent\hangindent
2\parindent\ignorespaces}

\long\def\eatit#1{}

\def\leftheadlinetext{}
\def\rightheadlinetext{}

\def\leftheadline{{\eightrm\folio\hfil \leftheadlinetext\hfil}}
\def\rightheadline{{\eightrm\hfil\rightheadlinetext\hfil\folio}}

\headline={\ifnum\pageno=1\hfil\else
\ifodd\pageno\rightheadline\else\leftheadline\fi\fi}

\def\tenpoint{\def\rm{\fam0\tenrm}
\textfont0=\tenrm \scriptfont0=\sevenrm \scriptscriptfont0=\fiverm
\textfont1=\teni \scriptfont1=\seveni \scriptscriptfont1=\fivei
\def\mit{\fam1} \def\oldstyle{\fam1\teni}
\textfont2=\tensy \scriptfont2=\sevensy \scriptscriptfont2=\fivesy
\def\cal{\fam2}
\textfont3=\tenex \scriptfont3=\tenex \scriptscriptfont3=\tenex
\def\it{\fam\itfam\tenit} 
\textfont\itfam=\tenit
\def\sl{\fam\slfam\tensl} 
\textfont\slfam=\tensl
\def\bf{\fam\bffam\tenbf} 
\textfont\bffam=\tenbf \scriptfont\bffam=\sevenbf
\scriptscriptfont\bffam=\fivebf
\def\tt{\fam\ttfam\tentt} 
\textfont\ttfam=\tentt
\normalbaselineskip=12pt
\setbox\strutbox=\hbox{\vrule height8.5pt depth3.5pt  width0pt}%
\normalbaselines\rm}

\def\eightpoint{\def\rm{\fam0\eightrm}%
\textfont0=\eightrm \scriptfont0=\sixrm \scriptscriptfont0=\fiverm
\textfont1=\eighti \scriptfont1=\sixi \scriptscriptfont1=\fivei
\def\mit{\fam1} \def\oldstyle{\fam1\eighti}%
\textfont2=\eightsy \scriptfont2=\sixsy \scriptscriptfont2=\fivesy
\def\cal{\fam2}%
\textfont3=\tenex \scriptfont3=\tenex \scriptscriptfont3=\tenex
\def\it{\fam\itfam\eightit} 
\textfont\itfam=\eightit
\def\sl{\fam\slfam\eightsl} 
\textfont\slfam=\eightsl
\def\bf{\fam\bffam\eightbf} 
\textfont\bffam=\eightbf \scriptfont\bffam=\sixbf
\scriptscriptfont\bffam=\fivebf
\def\tt{\fam\ttfam\eighttt} 
\textfont\ttfam=\eighttt
\normalbaselineskip=9pt%
\setbox\strutbox=\hbox{\vrule height7pt depth2pt  width0pt}%
\normalbaselines\rm}

\def\largefont{\def\rm{\fam0\largerm}
\textfont0=\largerm \scriptfont0=\largescriptrm \scriptscriptfont0=\tenrm
\textfont1=\largei \scriptfont1=\largescripti \scriptscriptfont1=\teni
\def\mit{\fam1} \def\oldstyle{\fam1\teni}
\textfont2=\largesy 
\def\cal{\fam2}
\def\it{\fam\itfam\largeit} 
\textfont\itfam=\largeit
\def\sl{\fam\slfam\largesl} 
\textfont\slfam=\largesl
\def\bf{\fam\bffam\largebf} 
\textfont\bffam=\largebf 
\scriptscriptfont\bffam=\fivebf
\def\tt{\fam\ttfam\largett} 
\textfont\ttfam=\largett
\normalbaselineskip=17.28pt
\setbox\strutbox=\hbox{\vrule height12.25pt depth5pt  width0pt}%
\normalbaselines\rm}

\def\titlefont{\def\rm{\fam0\titlerm}
\textfont0=\titlerm \scriptfont0=\largescriptrm \scriptscriptfont0=\tenrm
\textfont1=\titlei \scriptfont1=\largescripti \scriptscriptfont1=\teni
\def\mit{\fam1} \def\oldstyle{\fam1\teni}
\textfont2=\titlesy 
\def\cal{\fam2}
\def\it{\fam\itfam\titleit} 
\textfont\itfam=\titleit
\def\sl{\fam\slfam\titlesl} 
\textfont\slfam=\titlesl
\def\bf{\fam\bffam\titlebf} 
\textfont\bffam=\titlebf 
\scriptscriptfont\bffam=\fivebf
\def\tt{\fam\ttfam\titlett} 
\textfont\ttfam=\titlett
\normalbaselineskip=24.8832pt
\setbox\strutbox=\hbox{\vrule height12.25pt depth5pt  width0pt}%
\normalbaselines\rm}

\nopagenumbers

\font\eightrm=cmr8
\font\eighti=cmmi8
\font\eightsy=cmsy8
\font\eightbf=cmbx8
\font\eighttt=cmtt8
\font\eightit=cmti8
\font\eightsl=cmsl8
\font\sixrm=cmr6
\font\sixi=cmmi6
\font\sixsy=cmsy6
\font\sixbf=cmbx6

\font\largerm=cmr12 at 17.28pt
\font\largei=cmmi12 at 17.28pt
\font\largescriptrm=cmr12 at 14.4pt
\font\largescripti=cmmi12 at 14.4pt
\font\largesy=cmsy10 at 17.28pt
\font\largebf=cmbx12 at 17.28pt
\font\largett=cmtt12 at 17.28pt
\font\largeit=cmti12 at 17.28pt
\font\largesl=cmsl12 at 17.28pt

\font\titlerm=cmr12 at 24.8832pt
\font\titlei=cmmi12 at 24.8832pt
\font\titlesy=cmsy10 at 24.8832pt
\font\titlebf=cmbx12 at 24.8832pt
\font\titlett=cmtt12 at 24.8832pt
\font\titleit=cmti12 at 24.8832pt
\font\titlesl=cmsl12 at 24.8832pt

\tenpoint



\def\manyby{\hbox to.75in{\hrulefill}}
\hsize 5.41667in 
\vsize 7.5in

\def\manyby{\hbox to.75in{\hrulefill}}

\tolerance 3000
\hbadness 3000

\def\item#1{\par\indent\indent\llap{\rlap{#1}\indent}\hangindent
2\parindent\ignorespaces}

\def\refAHb{AH}
\def\refCat{Cat}
\def\refCMa{CM1}
\def\refCMb{CM2}
\def\refE{Ev}
\def\refGGR{GGR}
\def\refvanc{H1}
\def\refigp{H2}
\def\refHHF{HHF}
\def\refHR{HR}
\def\refHi{Hi1}
\def\refHib{Hi2}
\def\refId{Id}
\def\refNtwo{N}

\def\binom#1#2{\hbox{$#1 \choose #2$}}

\def\pr#1{\hbox{{\bf P}${}^{#1}$}}
\def\cite#1{[\ir{#1}]}

\def\leftheadlinetext{B. Harbourne and J. Ro\'e}
\def\rightheadlinetext{Fat Points on \pr2}

\title{Hilbert functions and resolutions
for ideals of $n=s^2$ fat points in \pr2}

\author{Brian Harbourne}

\address
Department of Mathematics and Statistics
University of Nebraska-Lincoln
Lincoln, NE 68588-0323
email: bharbour@math.unl.edu
WEB: {\tt http://www.math.unl.edu/$\sim$bharbour/}
\smallskip
\endaddress

\author{Joaquim Ro\'e}

\address
Departament d'\`Algebra i Geometria
Universitat de Barcelona
Barcelona 08007, Spain
email: jroevell@mat.ub.es
\smallskip

April 24, 2001\endaddress
\vskip-\baselineskip

\thanks{\vskip -6pt
\noindent 1980 {\it Mathematics Subject Classification. } 
Primary 13P10, 14C99. 
Secondary 13D02, 13H15.
\smallskip
\noindent {\it Key words and phrases. }  Hilbert function,
resolution, fat points, blow up, rational surface.\smallskip}

\vskip\baselineskip
\Abstract{Abstract: Conjectures for the Hilbert function
of the $m$th symbolic power of the ideal of $n$ 
general points of \pr2 are verified for infinitely many $m$ 
for each square $n>9$, using an approach developed by the authors 
in a previous paper. In those cases that $n$ is even, 
conjectures for the resolution are also verified. Previously,
by work of Evain (for the Hilbert function)
and of Harbourne, Holay and Fitchett (for the resolution),
these conjectures were known for infinitely many $m$ for a given
$n$ only for $n$ being a power of 4.}
\vskip\baselineskip

\irrnSection{Introduction}{intro}
Consider the ideal $I_{(n,m)}\subset R={\bf C}[\pr2]$ generated by 
all forms having multiplicity at least $m$ at $n$ given general points
of \pr2. This is a graded ideal, and thus we can consider the 
Hilbert function $h_{(n,m)}$ whose value at each nonnegative integer 
$t$ is the dimension $h_{(n,m)}(t) = \hbox{dim} (I_{(n,m)})_t$
of the homogeneous component $(I_{(n,m)})_t$ 
of $I_{(n,m)}$ of degree $t$.  

In spite of a good deal of work, until recently $h_{(n,m)}$
was known only when either $n$ or $m$ was small (see Figure 1,
a color postscript graphic, which can be viewed, in color,
also at http://www.math.unl.edu/$\sim$bharbour/Hilb1.jpg). Alexander 
and Hirschowitz \cite{refAHb} determined $h_{(n,m)}$ for each $m$
for all $n$ sufficiently large compared with $m$ (in all dimensions,
not just for \pr2); however, it is unclear how large is large enough.

The first explicit determination when both $m$ and $n$ can be large 
was given by Evain \cite{refE}, whose result applies for all $m$
whenever $n$ is a power of 4. (Evain's method also seems to work
as long as $n$ is a square divisible only by 2, 3 or 5.)
We \cite{refHR} recently found a different method and applied it
to determine $h_{(n,m)}$ in a great many additional cases
in which both $m$ and $n$ can be large, but, unlike Evain's result,
for these cases the larger $m$ is, the larger $n$ must be also.

Here we show that our method also applies, like Evain's, to determine
$h_{(n,m)}$ for infinitely many $m$ for various $n$ (in fact,
whenever $n$ is a square). In addition, applying the results of
\cite{refHHF}, we also obtain the resolution of the ideal $I_{(n,m)}$
for infinitely many $m$ whenever $n$ is an even square.
Thus with this paper and with our recent results, we can update the 
data of Figure 1; the updated data is shown in Figure 2
(on the web this is located at 
http://www.math.unl.edu/$\sim$bharbour/Hilb2.jpg). See Figures 3 
and 4 for the corresponding data for resolutions (on the web at 
http://www.math.unl.edu/$\sim$bharbour/Res1.jpg and
http://www.math.unl.edu/$\sim$bharbour/Res2.jpg).
[Note that these figures do not show all cases for which the 
algorithm of \cite{refHR} determines the Hilbert function or resolution,
just certain cases for which the algorithm is especially easy 
to analyze (i.e., Corollaries V.2 and V.4 of \cite{refHR}) and those 
additional cases we analyze here (i.e., \ir{mainthm}). 
There are additional cases
which a computer search shows that the algorithm handles, but
for which a simple statement seems hard to give.]

\irrnSection{Conjectures and Main Result}{main}
As shown in Figure 1, $h_{(n,m)}$ was determined for all $m$
by Nagata, when $n<10$. The result in these cases turns out to be
a bit complicated. Conjectures have been given
(\cite{refvanc},\cite{refHi}) which imply
that these complications disappear for $n\ge10$:

\prclm{Conjecture}{Hconj}{For $n\ge10$ general points of \pr2, 
$h_{(n,m)}(t)=\hbox{max}\{0, \binom{t+2}{2}-n\binom{m+1}{2}\}$
for each integer $t\ge0$.}

Similarly, as shown in Figure 3, the resolution of $I_{(n,m)}$ is known
(\cite{refCat}, \cite{refigp}) when $n\le9$, but it also is somewhat 
complicated. These complications are conjectured to disappear for
$n\ge10$ (\cite{refigp}, \cite{refHHF}):

\prclm{Conjecture}{Rconj}{For $n\ge10$, the minimal free 
resolution of $I_{(n,m)}$ is 
$$0\to R[-\alpha-2]^d\oplus R[-\alpha-1]^c\to 
R[-\alpha-1]^b\oplus R[-\alpha]^a\to I_{(n,m)}\to 0,$$
where $\alpha$ is the least $t$ such that $h_{(n,m)}(t)>0$,
$a=h_{(n,m)}(\alpha)$, 
$b=\hbox{max}\{h_{(n,m)}(\alpha+1)-3h_{(n,m)}(\alpha),0\}$,
$c=\hbox{max}\{-h_{(n,m)}(\alpha+1)+3h_{(n,m)}(\alpha),0\}$,
$d=a+b-c-1$, and $R[i]^j$ is the direct sum of $j$ copies of
the ring $R={\bf C}[\pr2]$, regarded as an $R$-module
with the grading $R[i]_k=R_{k-i}$.}

We now show that \ir{Hconj} holds for infinitely
many $m$ whenever $n\ge10$ is a square, and (as 
a consequence of \cite{refHHF}) that \ir{Rconj} holds for infinitely many
$m$ whenever $n\ge10$ is an even square. For the purpose of 
stating the theorem, given any positive integer $i$, let
$l_i$ be the largest integer $j$ such that $j(j+1)\le i$.

\prclm{Theorem}{mainthm}{Consider $10\le n=s^2$ general points of \pr2.
Let $k$ be any nonnegative integer, and let $m=x+k(s-1)$,
where $x$ is an integer satisfying $s/2-l_{s}\le x\le s/2$
if $s$ is even, or $(s+1)/2-l_{2s}\le x\le (s+1)/2$
if $s$ is odd. Then \ir{Hconj} holds for $I_{(n,m)}$ if $n$ is odd,
and \ir{Hconj} and \ir{Rconj} both hold if $n$ is even.}

To prove \ir{mainthm} we will, in the particular case of 
uniform multiplicities, use an algorithm developed in
\cite{refHR}. The algorithm gives bounds on 
$\alpha(n,m)$, defined to be the least $i$ such that 
$h_{(n,m)}(i)>0$. We begin by briefly
recalling the algorithm. The parenthetical remarks below
can be ignored; they are meant only as a reminder of the 
geometry underlying the algorithm as developed in \cite{refHR}, 
and are not necessary for applying the algorithm.

Let $g=(d-1)(d-2)/2$ (thus $g$ is the genus of a 
plane curve of degree $d$).
Let $C$ denote the $n+1$-tuple $(d,1,\ldots,1,0,\ldots,0)$, where
the entry 1 occurs $r$ times. Given any $n+1$-tuple
$D=(t,m_1,m_2,\ldots, m_n)$ of integers,
with respect to the obvious bilinear form we have
$D\cdot C=td-(m_1+\cdots+m_r)$. (The $n+1$-tuple
$D$ corresponds to a divisor ${\tilde D}$ on
the blow up $X$ of \pr2 at the $n$ general points,
$C$ corresponds to the proper transform ${\tilde C}$ 
of a curve of
degree $d$ passing through $r$ of the points and
$D\cdot C$ is just the usual intersection product
${\tilde D}\cdot{\tilde C}$, and hence equals
the degree of ${\cal O}_{\tilde C}({\tilde D})$.)

Assume that $D_0=(t,m_1,m_2,\ldots, m_n)$ satisfies
$m_1\ge m_2\ge \cdots\ge m_n\ge0$.
Define $D_i'$ by $D_i'=D_{i-1}-C$, and then
define $D_i$ to be what is obtained from $D_i'$ by first
converting any of the 2nd through $n+1$st entries to 0 if
it is negative and then permuting the 2nd through $n+1$st
entries to be in descending order. For convenience, define
$t_i$ to be the first entry of $D_i$ (thus 
$t_i={\tilde D_\omega}\cdot L$ where $L$ is the total transform
to $X$ of a line in \pr2).

We thus get a sequence of tuples, $D_0$, $D_1$, etc.;
let $\omega$ be the least $i$ such that  $t_i<0$. By \cite{refHR},
we then have $\alpha(n,m)>t$ if the following two conditions
hold: (1) $D_i\cdot C\le g-1$ for $0\le i<\omega-1$, and
(2) $(t-({\omega-1})d+1)(t-({\omega-1})d+2)\le 2\mu$,
where $\mu=t_{\omega-1}d-D_{\omega-1}\cdot C$ is 
the sum of the second to $r$th 
entries of $D_{\omega-1}$. 

The algorithm is simply to determine the 
biggest $t$ such that the given conditions are satisfied.
(Geometrically, $t_\omega={\tilde D_\omega}\cdot L<0$
and hence ${\tilde D_\omega}$ is not
linearly equivalent to an effective divisor.
According to the analysis in \cite{refHR},
$D_i\cdot C\le g-1$ for $i<\omega-1$ and
$(t_i+1)(t_i+2)\le 2\mu$ for $i=\omega-1$
mean that ${\cal O}_{\tilde C}({\tilde D}_i)$
has no global sections,
which by induction implies the divisor
corresponding to ${\tilde D}_0$ is not
linearly equivalent to an effective divisor.)

In order to analyze the algorithm, we will use 
the following lemma, for which we define $\omega'$
to be the least $i$ such that all entries of $D_i$ except 
possibly the first are 0.

\prclm{Lemma}{intlemma}{Let $D_0=(t,m,m,\ldots, m)$. Then
for all $0\le i\le \omega'-1$
$$
i\left({{r^2}\over{n}}-d^2\right)-\left(r-{{r^2}\over{n}}\right)
\le D_i \cdot C - D_0\cdot C \le
i\left({{r^2}\over{n}}-d^2\right).
$$
}

\noindent {\bf Proof}:
Let $A_0$ be the $n+1$-tuple $(0,\ldots,0)$,
and for $0< k\le n$ let
$A_k=(0,1,\ldots,1,0,\ldots,0)$, where 1's occupy
positions $1$ through $k$. 
For $0\le i<\omega'$, it is not hard to check that
$D_i=(t_0-id,m-i+q,\ldots,m-i+q)+A_\rho$,
where $i(n-r)=qn+\rho$ with $0\le\rho<n$,
and therefore
$$
D_i \cdot C - D_0 \cdot C= i(r-d^2) - rq +A_\rho \cdot C.
$$
On the other hand, $A_\rho \cdot C=-\min\{\rho,r\}$,
and it is easy to see that
$$
i(n-r){r\over n} \le rq+\min\{\rho,r\} \le
i(n-r){r\over n}+r-{{r^2}\over n}
$$
from which the claim follows.
\qed

Our next result improves (for uniform multiplicities)
on Theorem I.1(a)(i) of \cite{refHR},
which assumed that $r \le d^2$, whereas here we only require that
$r^2 \le d^2 n$.

\prclm{Proposition}{aprop}{Consider $n$ general points of \pr2.
Let $m$, $r$, $d$, $u$ and $\rho$ be nonnegative integers such that
$nm=ur+\rho$, $0<\rho\le r\le n$, and $rd(d+1)/2\le r^2 \le d^2 n$,
and let $l=\min\{l_{2\rho},d\}-1$. Then $\alpha{(n,m)}\ge 
1+\min \{\lfloor (mr+g-1)/d \rfloor, l+ud\}$.}

\noindent {\bf Proof}: It is easy to check that 
$\omega'=\lceil mn/r\rceil=u+1$,
so if $t \le l+ud$, it follows that $t_{\omega'}=l-d<0$,
and thus $\omega'\ge\omega$. Since $r^2/n-d^2 \le 0$, 
it follows from \ir{intlemma} that
$D_i\cdot C\le D_0 \cdot C$ for all $0\le i\le \omega-2$.
If $t\le \lfloor (mr+g-1)/d \rfloor$, then
$D_i\cdot C\le D_0 \cdot C = td-mr \le g-1$.
To conclude that $\alpha(n,m)\ge t+1$, it is now enough to 
check that $(t-i'd+1)(t-i'd+2)\le 2\mu$ for $i'=\omega-1$,
where $\mu$ is the sum of the second to $r$th entries of $D_{i'}$.
If $i'=u$ (i.e., $\omega'=\omega$)
we have $t-i'd=t-ud\le l$ by hypothesis and hence 
$(t-i'd+1)(t-i'd+2)\le (l+1)(l+2)\le 2\rho=2\mu$ by 
definition of $l$. If $i'<u$ (so  $\omega'>\omega$), 
by definition of $i'$ we at least have $t-i'd\le d-1$, so 
$(t-i'd+1)(t-i'd+2)\le d(d+1)$.
But $\omega'>\omega$ implies $\mu\ge r$, and
by hypothesis $rd(d+1)/2\le r^2$ (so $d(d+1)\le 2r$);
therefore $2\mu\ge 2r\ge d(d+1)\ge (t-i'd+1)(t-i'd+2)$
as we wanted. \qed

We now give the proof of our main result:

\noindent {\bf Proof of \ir{mainthm}}:
We apply \ir{aprop} with $d=s-1$, $r=ds$,
$u=\lceil mn/r\rceil-1=m+k$ and
$\rho=mn-ur=xs$. We claim
that $t_0 \le \min \{\lfloor (mr+g-1)/d \rfloor, l+ud\}$,
where $t_0=ms+s/2-2$ if $s$ is even 
and $t_0=ms+(s-1)/2-2$ if $s$ is odd.
But $t_0\le(mr+g-1)/d$ because
$(mr+g-1)/d=(mds+d(d-3)/2)/d=ms+(s-4)/2$.
To see $t_0 \le l+ud$, note that $t_0 \le l+ud$ 
simplifies to $x+s/2-2 \le l$ if $s$ is even and to
$x+(s-1)/2-2 \le l$ if $s$ is odd.
Therefore (by definition of $l$) we have to
check that $x+s/2-1 \le d$ and $(x+s/2-1)(x+s/2)\le 2sx$
if $s$ is even, and that $x+(s-1)/2-1 \le d$ and
$(x+(s-1)/2-1)(x+(s-1)/2)\le 2sx$ if $s$ is odd.
The first inequality follows from $x\le s/2$ and
$x \le (s+1)/2$ respectively. For the second,
substituting $s/2-j$ for $x$ if $s$ is even and
$(s+1)/2-j$ for $x$ if $s$ is odd,
$(x+s/2-1)(x+s/2)\le 2sx$ and $(x+(s-1)/2-1)(x+(s-1)/2)\le 2sx$ resp.
become $j(j+1)\le s$ if $s$ is even and $j(j+1)\le 2s$ if
$s$ is odd. Thus $(x+s/2-1)(x+s/2)\le 2sx$ and 
$(x+(s-1)/2-1)(x+(s-1)/2)\le 2sx$ resp.
hold if $x$ is an integer satisfying $s/2-l_{s}\le x\le s/2$
if $s$ is even, and $(s+1)/2-l_{2s}\le x\le (s+1)/2$
if $s$ is odd.

This shows by \ir{aprop} that $\alpha(n,m)\ge ms+s/2-1$ 
if $s$ is even and
$\alpha(n,m)\ge ms+(s-1)/2-1$ if $s$ is odd. But
since $n$ points of multiplicity $m$ impose at most
$n\binom{m+1}{2}$ conditions on forms of degree $t$, it
follows that $h_{(n,m)}(t)\ge \binom{t+2}{2}-n\binom{m+1}{2}$,
and it is easy to check  that $\binom{t+2}{2}-n\binom{m+1}{2}>0$ 
whenever $t\ge ms+s/2-1$ if $s$ is even 
and $t\ge ms+(s-1)/2-1$ if $s$ is odd.
Thus in fact we have $\alpha(n,m)=ms+s/2-1$ if $s$ is even and
$\alpha(n,m)=ms+(s-1)/2-1$ if $s$ is odd, whenever $m$ is of the
form $m=x+k(s-1)$, with $x$ as given in the statement
of \ir{mainthm}.

Of course, $h_{(n,m)}(t)=0$ for all $t<\alpha(n,m)$, and
by \cite{refHHF}, we know that $h_{(n,m)}(t)=
\binom{t+2}{2}-n\binom{m+1}{2}$
for all $t\ge\alpha(n,m)$ (apply Lemma 5.3 of \cite{refHHF},
keeping in mind our explicit expression for $\alpha(n,m)$). 
This proves our claims regarding
verification of \ir{Hconj}. As for \ir{Rconj}, when $n$
is an even square, apply Theorem 5.1(a) of \cite{refHHF}.
This concludes the proof. \qed

\irrnSection{Figures}{figs}
Here we show graphs of what was known until recently regarding
\ir{Hconj} and \ir{Rconj}, and what is known now. For Figure 1,
the references are: \cite{refNtwo} for Nagata, 
\cite{refHib} for Hirschowitz, \cite{refCMa} and \cite{refCMb}
for  Ciliberto and Miranda, and \cite{refE} for Evain.
The additional data shown in Figure 2 is simply a graphical 
representation of Corollary V.2 of \cite{refHR},
and \ir{mainthm}.

For Figure 3, the references are: \cite{refGGR} for Geramita, 
Gregory and Roberts, \cite{refCat} for Catalisano, \cite{refigp} for 
Harbourne, \cite{refId} for Id\`a, and \cite{refHHF} for Harbourne, 
Holay, and Fitchett. As before, 
the additional data shown in Figure 4 is simply a graphical 
representation of Corollary V.4 of \cite{refHR},
and \ir{mainthm}. 

\vfill\eject
\noindent\hbox to\hsize{\hfil Figure 1\hfil}
\vskip7in
\includegraphics{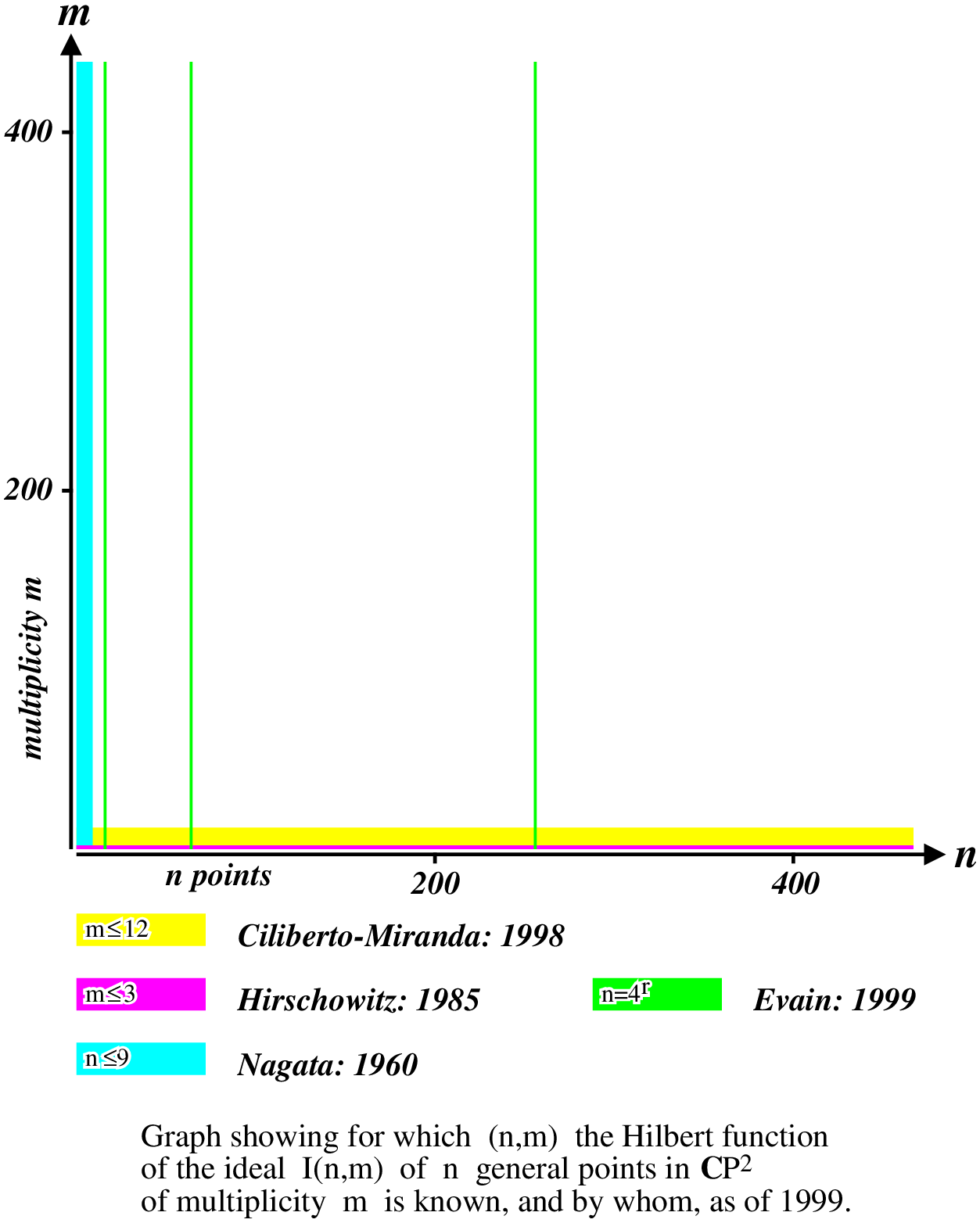}
\vfill\eject
\noindent\hbox to\hsize{\hfil Figure 2\hfil}
\vskip\vsize
\vskip-\baselineskip
\includegraphics{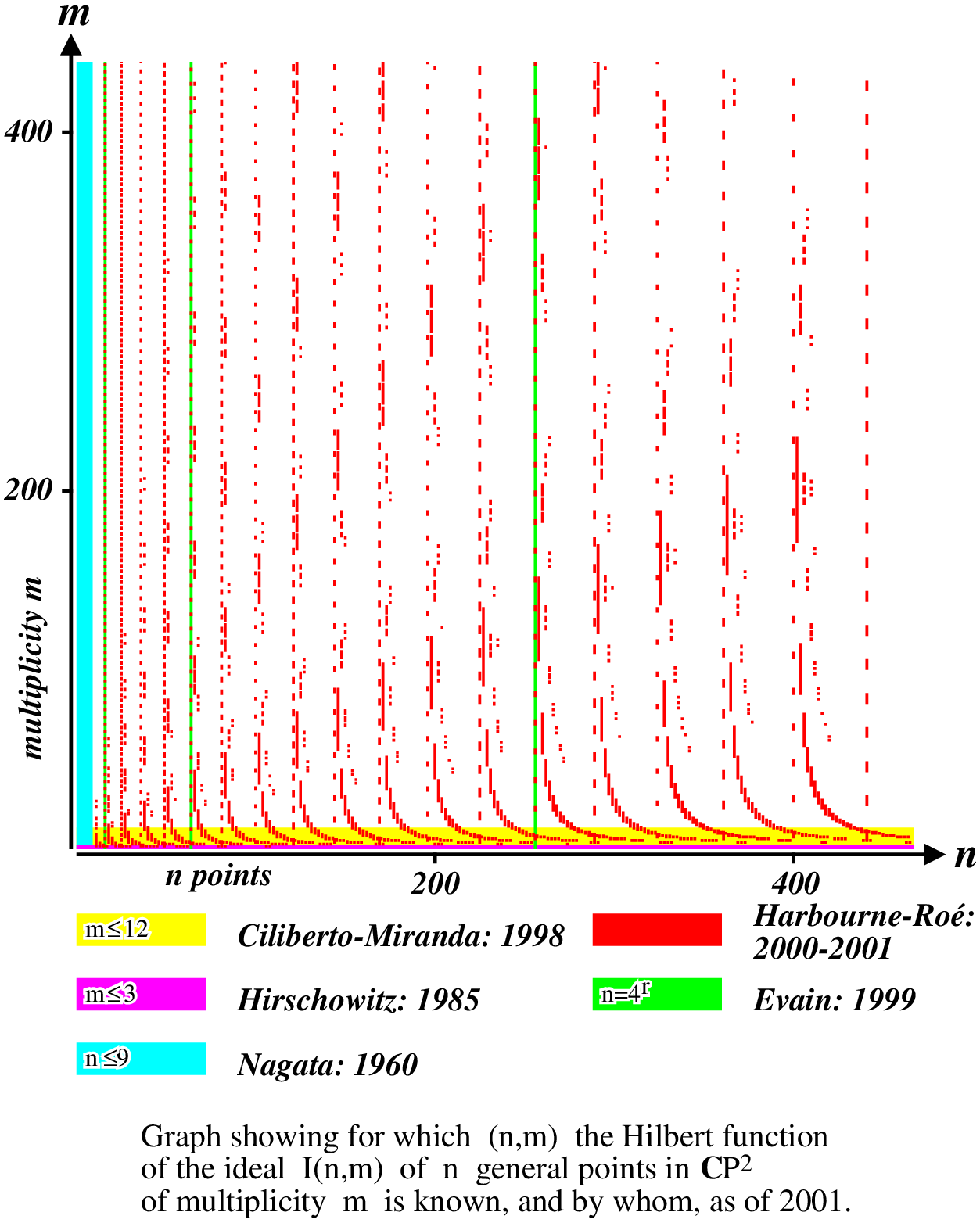}
\vfill\eject
\noindent\hbox to\hsize{\hfil Figure 3\hfil}
\vskip\vsize
\vskip-\baselineskip
\includegraphics{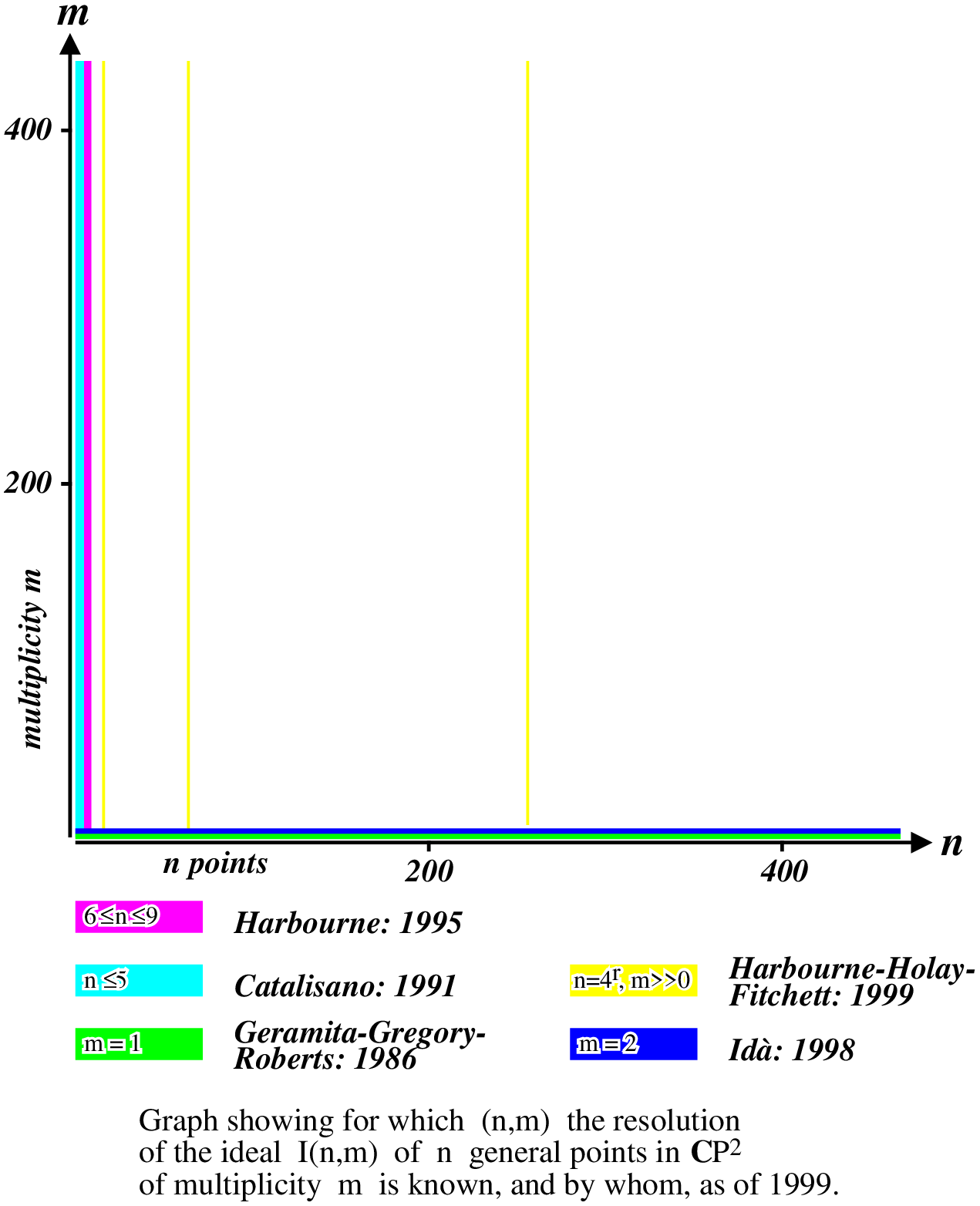}
\vfill\eject
\noindent\hbox to\hsize{\hfil Figure 4\hfil}
\vskip\vsize
\vskip-\baselineskip
\includegraphics{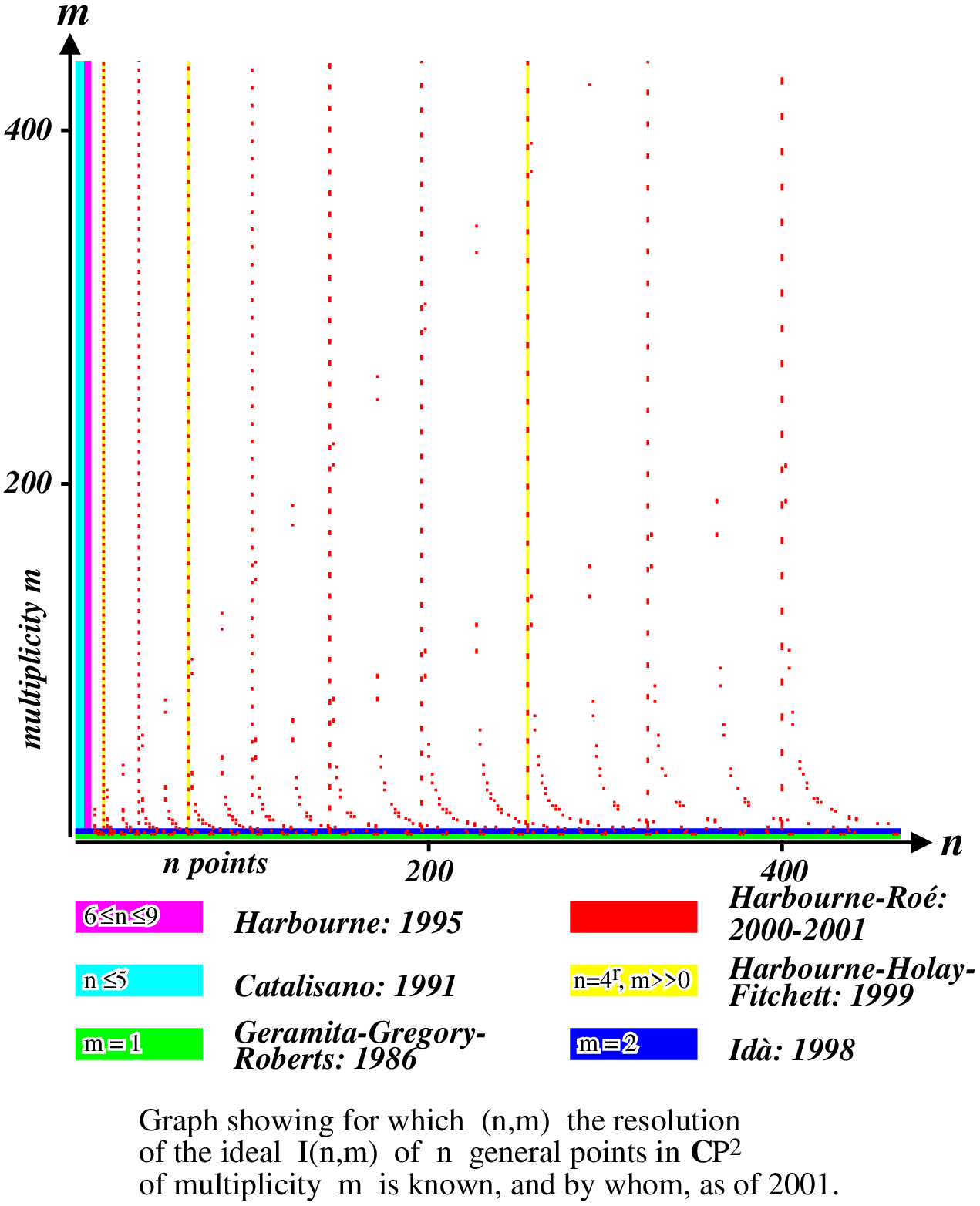}
\vfill\eject

\References

\bibitem{\refAHb} J. Alexander and A. Hirschowitz. {\it An asymptotic 
vanishing theorem for generic unions of multiple points}, 
Invent. Math. 140 (2000), no. 2, 303--325.

\bibitem{\refCMa} C. Ciliberto and R. Miranda. {\it Degenerations
of planar linear systems}, 
J. Reine Angew. Math. 501 (1998), 191-220.

\bibitem{\refCMb} C. Ciliberto and R. Miranda. {\it Linear systems
of plane curves with base points of equal multiplicity}, 
Trans. Amer. Math. Soc. 352 (2000), 4037--4050.

\bibitem{\refCat} M.\ V. Catalisano. {\it ``Fat'' points on a conic}, 
Comm.\ Alg.\  19(8) (1991), 2153--2168.

\bibitem{\refE} L. Evain.
{\it La fonction de Hilbert de la r\'eunion
de $4^h$ gros points g\'en\'eriques
de \pr2 de m\^eme multiplicit\'e},
J. Alg. Geom. 8 (1999), 787--796.

\bibitem{\refGGR} A.\ V. Geramita, D. Gregory and L. Roberts.
{\it Monomial ideals and points in projective space},
J.\ Pure and Appl.\ Alg. 40 (1986), 33--62.

\bibitem{\refvanc} B. Harbourne. {\it The geometry of 
rational surfaces and Hilbert
functions of points in the plane},
Can.\ Math.\ Soc.\ Conf.\ Proc.\ 6 
(1986), 95--111.

\bibitem{\refigp} \manyby. {\it  The Ideal Generation 
Problem for Fat Points},
J. Pure and Applied Alg. 145 (2000), 165--182.

\bibitem{\refHHF} B. Harbourne, S. Holay, and S. Fitchett.
{\it Resolutions of Ideals of Quasiuniform Fat Point Subschemes of \pr2},
preprint (2000).

\bibitem{\refHR} B. Harbourne and J. Ro\'e. 
{\it Linear systems with multiple base points in \pr2}, preprint (2000).

\bibitem{\refHi} A. Hirschowitz.
{\it Une conjecture pour la cohomologie 
des diviseurs sur les surfaces rationelles g\'en\'eriques},
Journ.\ Reine Angew.\ Math. 397
(1989), 208--213.

\bibitem{\refHib} A. Hirschowitz.
{\it La m\'ethode d'Horace pour l'interpolation \`a plusieurs
variables}, Manus. Math. 50 (1985), 337--388.

\bibitem{\refId} M. Id\`a.
{\it The minimal free resolution for the first infinitesimal 
neighborhoods of $n$ general points in the plane},
J.\ Alg. 216 (1999), 741--753.

\bibitem{\refNtwo} M. Nagata. {\it On rational surfaces, II}, 
Mem.\ Coll.\ Sci.\ Univ.\ Kyoto, Ser.\ A Math.\ 33 (1960), 271--293.

\bye